\documentclass[12pt]{report}
\usepackage{amssymb,amsmath,makeidx,verbatim}
\newcommand{\R}{\mathbb{R}}

\newcommand{\be}{\begin{enumerate}}
	\newcommand{\ee}{\end{enumerate}}
\newcommand{\bq}{\begin{eqnarray*}}
	\newcommand{\eq}{\end{eqnarray*}}

\begin{document}
	\newcommand{\disp}{\displaystyle}
	\thispagestyle{empty}
	\begin{center}
		\textsc{WHAT IS.........\\
			........ the JEFT?\\
		\ \\
		\ \\
		\ \\
	\ \\
by}
		\ \\
		\ \\
		\textsc{Olufemi O. Oyadare}\\
		\ \\
		Department of Mathematics,\\
		Obafemi Awolowo University,\\
		Ile-Ife, $220005,$ NIGERIA.\\
		\text{E-mail: \textit{femi\_oya@yahoo.com}}\\
	\end{center}
		\ \\
		\ \\
	\ \\
	\indent The {\it JEFT} is the {\it acronym} for the Joint-Eigenspace Fourier Transform defined on a noncompact symmetric space. It is a consequence of a general construction of a Fourier transform modelled on the Harish-Chandra Fourier transform (on a semi-simple Lie group with finite centre) which (on the corresponding symmetric space of the noncompact type) serves as the Poisson-completion of the famous Helgason Fourier transform.
	
	\indent For a noncompact semi-simple Lie group $G$ (with finite centre) whose corresponding Lie algebra $\mathfrak{g}$ has the Cartan decomposition $\mathfrak{g}=\mathfrak{t}\oplus\mathfrak{p},$ its Iwasawa decomposition is given as $G=KAN$ in which $K$ is the analytic subgroup of $G$ with Lie algebra $\mathfrak{t},$ $A=:\exp(\mathfrak{a})$ (where $\mathfrak{a}$ is a maximal abelian subspace of $\mathfrak{p}$) and $N$ is the analytic subgroup of $G$ corresponding to $\mathfrak{n}=\sum_{\lambda\in\triangle^{+}}\{X\in\mathfrak{g}:[H,X]=\lambda(H)X,\forall H\in \mathfrak{a}\}$ (where $\triangle^{+}$ denote the set of all restricted positive roots). A member $\varphi\in C(G),$ with $\varphi(e)=1$ in which $\varphi(k_{1}gk_{2})=\varphi(g),$ for all $k_{1},k_{2}\in K, g\in G$ is termed a spherical function and is said to belong to $C(G//K).$ The Harish-Chandra spherical transform (on $C(G//K)$) written as $f\mapsto\widehat{f}$ is defined on $\mathfrak{a}^{*}_{\mathbb{C}}$ as $\widehat{f}(\lambda)=(f*\varphi_{\lambda})(e).$ Here $*$ is the convolution on $G,$ $\varphi_{\lambda}$ is the elementary spherical function corresponding to $\lambda$ and $e$ is the identity of $G.$ It would be more satisfying to consider a general (Fourier) transform $f\mapsto\mathcal{H}_{g}f$ of $C(G)$ defined on $\mathfrak{a}^{*}_{\mathbb{C}}$ as $(f*\varphi_{\lambda})(g),$ for every $g\in G$ (not just for $g=e$ of the above Harish-Chandra case) and to study the extent to which it could generalize the spherical transform on $G.$ We believe that every Fourier transform on any mathematical object should be sought in the above closed-convolution form.
	
	\indent Notable among the properties of the spherical transform is the Trombi-Varadarajan theorem which gives a characterization of the spherical transform image of the Schwartz Frechet algebra $\mathcal{C}^{p}(G//K)$ for all $0<p\leq2, [8.].$ The complete result $\mathcal{C}^{2}(SL(2,\mathbb{R}))$ and $\mathcal{C}^{p}(SL(2,\mathbb{R}))$ are contained in $[9.]$ and $[1.].$
	
	\indent Intimate relationship between semi-simple Lie groups (on one hand) and symmetric spaces in the sense of E. Cartan (on the other) suggests an equivalent theory of Fourier transform on the symmetric spaces which are found to be the quotient space $G/K=:X,$ for various pairs $(G,K).$ Indeed, if $M=\{k\in K: Ad(k)H=H,\;H\in \mathfrak{a}\}$ denote the centralizer of $\mathfrak{a}$ in $K,$ then the Helgason Fourier transform $f\mapsto\widetilde{f}$ (for every $f\in C^{\infty}_{c}(X)$) is defined on $\mathfrak{a}^{*}_{\mathbb{C}}\times B$ (with $B:=K/M$) as the integral $$\widetilde{f}(\lambda,b)=\int_{X}f(x)e^{(-i\lambda+\rho)(A(x,b))}dx,$$ whenever it exists. Here $A(x,b)=A(gK,kM)=A(k^{-1}g),$ where $A(g):=-H(g^{-1})$ in which $H(g)\in \mathfrak{a}$ is determined by the Iwasawa decomposition $g=k\exp H(g)\mathfrak{n}$ with $\rho=\frac{1}{2}\sum_{\lambda\in\triangle^{+}}m_{\lambda}\cdot\lambda$ (in which $$m_{\lambda}=dim(\{X\in\mathfrak{g}:[H,X]=\lambda(H)X,\forall H\in \mathfrak{a}\}).$$
	
	\indent The Helgason Fourier transform enjoys much of the properties of the Harish-Chandra spherical transform, $[7.],$ notably as an isometry of $L^{2}(X)$ onto $L^{2}(\mathfrak{a}^{*}_{+}\times B)$ with Plancherel measure $|c(\lambda)|^{-2}d\lambda db,$ $[4.].$ The Helgason Fourier transform is indeed an extension (though not a generalization) of the Harish-Chandra spherical transform in the following sense.
	
	\indent {\bf Lemma 1} ($[3.],$ p. $224$) The Helgason Fourier transform is an extension of the Harish-Chandra spherical transform from $K-$invariant functions on $X$ to not-necessarily $K-$invariant functions on $X.\;\Box$
	
	\indent This lemma simply says that the Helgason Fourier transform is defined (whenever the above integral exists) for all functions in $C(X)=C(G/K)$ and not just for its subspace $C(G//K),$ but that whenever the Helgason Fourier transform is restricted to functions $f\in C(G//K)=C(K\backslash G/K)=C(K\backslash X),$ it coincides with the Harish-Chandra spherical transform of $f;$ directly implying that the restriction of the Helgason Fourier transform to $C(K\backslash X)$ has the closed-convolution form $(f\times\varphi_{\lambda})(eK),$ where the spherical function $\varphi_{\lambda}$ is now defined on $X$ as $\varphi_{\lambda}(gK)=\int_{K}e^{(i\lambda+\rho)(A(kg))}dk$ and $\times$ is the convolution on $X.$ It is however not valid to say that the general (i.e., unrestricted) Helgason Fourier transform has the closed-convolution form. Indeed, the general (unrestricted) Helgason Fourier transform on a noncompact symmetric space does not restrict to the general (unrestricted) classical Fourier transform on $\mathbb{R}^{n}.$
	
	\indent This may suggest that there is something inadequate or incomplete in the composition of the Helgason Fourier transform on $X,$ which needs to be fixed before the image $C^{\infty}_{c}(X)^{\sim}$ could be efficiently characterized. It is therefore appropriate to take a cue from the classical case of $X=\mathbb{R}^{n}$ in order to discuss more about the incompleteness of the Helgason Fourier transform.
	
	\indent The characters of $\mathbb{R}^{n}$ which double as the spherical functions are defined as $$\varphi_{t}:x\mapsto\varphi_{t}(x):=e^{-i<t,x>},\;t\in\mathbb{R}$$ in which $<t,x>=<(t_{1},\cdots,t_{n}),(x_{1},\cdots,x_{n})>:=t_{1}x_{1}+\cdots+t_{n}x_{n}.$ The Fourier transform of $f$ on $\mathbb{R}^{n}$ is therefore given as the convolution integral $\widetilde{f}(t)=\int_{\mathbb{R}^{n}}f(x)\varphi_{t}(x)dx,\;t\in\mathbb{R},$ and it maps the space $(L^{2}(\mathbb{R}),dx)$ onto $(L^{2}(\mathbb{R}^{n}),(2\pi)^{-n}dt)$ as a unitary isomorphism. For $\alpha=(\alpha_{1},\cdots\alpha_{n})\in\mathbb{N}\cup\{0\},$ set $D^{\alpha}:=(-i\frac{\partial}{\partial t_{1}})^{\alpha_{1}}\cdots(-i\frac{\partial}{\partial t_{n}})^{\alpha_{n}}$ and observe that $(D^{\alpha}\varphi_{t})(x)=\lambda_{\alpha}(x)\varphi_{t}(x),$ where $\lambda_{\alpha}(t)=\lambda_{\alpha}(t_{1},\cdots,t_{n}).$ The last differential equation implies that $(D^{\alpha}\widetilde{f}(t))(x)=\lambda_{\alpha}(x)\widetilde{f}(t).$ This informs that the classical Fourier transform on $\mathbb{R}^{n}$ is a joint-Eigenfunction, a property which is not {\it generally} attainable by the Helgason Fourier transform. These observations suggest that the Helgason Fourier transform needs a completion into becoming a joint-Eigenfunction.
	
	\indent The choice of the completion of the Helgason Fourier transform into the joint-Eigenspace is suggested by the following classical result on $\mathbb{R}^{n},$ which says that {\it the classical Fourier transform is the one-dimensional Fourier transform image of the Radon transform on $\mathbb{R}^{n}.$} That is, for any function $f$ on $\mathbb{R}^{n}$ that is integrable on each hyperplane, we have that $$\widetilde{f}(s\omega)=\int_{\mathbb{R}}\widehat{f}(\omega,p)e^{-ips}dp,$$ where $\widehat{f}(\omega,p)$ would be the Radon transform of $f,$ $[3.],$ p. $11.$ This would suggest that the classical Fourier transform $\widetilde{f}(s\omega)$ is a complete Fourier transform in that it is the one-dimensional-completion of the classical Radon transform. Thus the classical Radon transform is therefore not fully (or completely) a Fourier transform (until it is subjected to the one-dimensional completion given above), just as earlier observed about the Helgason Fourier transform. It is then imperative to ask: {\it what do the above one-dimensional Fourier transform and the classical Radon transform generalize into for a symmetric space of the noncompact type?} We shall thus quickly motivate the notion of the Radon transform and the one-dimensional Fourier transform in the general context of a symmetric space of the noncompact type, through which it would be clear that the Helgason Fourier transform on $X$ is more of an {\it offspring} of the classical Radon transform than being of the classical Fourier transform.
	
	\indent Without much ado, the Poisson transform $F\mapsto\mathcal{P}_{\lambda}F$ of $F\in C(B)$ is defined by $(\mathcal{P}_{\lambda}F)(x)=\int_{B}e^{(i\lambda+\rho)(A(x,b))}F(b)db,$ whenever this integral exists. It is known that $\mathcal{P}_{\lambda}$ maps $C(B)$ into the joint-eigenspace $$\mathcal{E}_{\lambda}(X):=\{f\in C^{\infty}(X): Df=\Gamma(D)(i\lambda)f\;\mbox{for }\;D\in \textbf{D}(X)\},$$ where $\textbf{D}(X)$ is the algebra of differential operators on $X$ that are invariant under all translations. Since we may consider the Helgason Fourier transform map as $f\mapsto\widetilde{f}_{\lambda}:b\mapsto\widetilde{f}_{\lambda}(b):=\widetilde{f}(\lambda,b),$ i.e., as a function of $b\in B,$ the Poisson transform image of the Helgason Fourier transform {\it may} be all that is needed to have {\it the} complete Fourier transform on any noncompact symmetric space $X.$ The foregoing is the motivation for the introduction of the JEFT. {\it What really is the JEFT?}
	
	\indent {\bf Definition}. Let $x\in X=G/K$ and let $f\in C^{\infty}_{c}(X).$ The transformation map $f\mapsto f^{\triangle}(\lambda,x)$ defined as $f^{\triangle}(\lambda,x):=(f\times \varphi_{\lambda})(x)$ is called the Joint-Eigenspace Fourier transform ({\it the JEFT})$.\;\Box$
	
	\indent As it is expected (to be a generalization of the Helgason Fourier transform, which is a function of $(\lambda,b)\in\mathfrak{a}_{\mathbb{C}}\times B\subset\mathfrak{a}_{\mathbb{C}}\times X$), the JEFT is a function of all of $(\lambda,x)\in \mathfrak{a}_{\mathbb{C}}\times X.$ This immediately reveals the self-duality of $X$ under the JEFT. Much more, the JEFT is in the closed-convolution form. The first pertinent questions are: {\it Is the JEFT indeed a Fourier transform on $X$ and is it distinct from the Helgason Fourier transform?}
	
	\indent An explication of the defining integral of the convolution $(f\times \varphi_{\lambda})(x)$ gives that $$f^{\triangle}(\lambda,x)=\int_{K/M}\int_{G/K}e^{(i\lambda+\rho)(A(x,b))}e^{(-i\lambda+\rho)(A(y,b))}f(y)dydb,$$ showing that $f\mapsto f^{\triangle}(\lambda,x)$ is a genuine Fourier transform on $C^{\infty}_{c}(X)$ into $C(\mathfrak{a}_{\mathbb{C}}\times X).$ The JEFT has the added advantage of being considered a function of $\lambda$ (as $f\mapsto f^{\triangle}_{x}:\lambda\mapsto f^{\triangle}_{x}(\lambda):=f^{\triangle}(\lambda,x)$) and (more importantly) as a function of $x$ (as $f\mapsto f^{\triangle}_{\lambda}:x\mapsto f^{\triangle}_{\lambda}(x):=f^{\triangle}(\lambda,x),$ on the entirety of $X$ as against Helgason Fourier transform's dependence on only $b\in B$). We may then use the above explicit form of the JEFT to study the Fourier transform theory for the compact symmetric space $K/M$ by restricting $f\mapsto f^{\triangle}(\lambda,x)$ to only $x=b\in K/M.$ What indeed is $f^{\triangle}(\lambda,eK)?$
	
	\indent Having shown that the JEFT is a genuine Fourier transform on $X,$ it is instructive to investigate how distinct it is from the Helgason Fourier transform (defined on the same $X$). {\it Or can we afford to have {\it two} distinct Fourier transforms on $X?$}
	
	\indent The definite answer to these questions are gleaned from the following Fubini argument, $[2.],$ p. $98.$
	
	\indent $$f^{\triangle}(\lambda,x)=\int_{K/M}\int_{G/K}e^{(i\lambda+\rho)(A(x,b))}e^{(-i\lambda+\rho)(A(y,b))}f(y)dydb,$$ $$\;\;\;\;\;\;\;\;\;\;\;\;\;\;\;
	=\int_{K/M}e^{(i\lambda+\rho)(A(x,b))}(\int_{G/K}e^{(-i\lambda+\rho)(A(y,b))}f(y)dy)db$$ $$\;\;\;\;\;\;\;\;\;\;\;
	=\int_{K/M}e^{(i\lambda+\rho)(A(x,b))}\widetilde{f}(\lambda,b)db=\mathcal{P}_{\lambda}(\widetilde{f}(\lambda,\cdot))(x).$$ We thus have the following.
	
	\indent {\bf Lemma 2} ($[5.],$ p. $12$). $f^{\triangle}(\lambda,x)=\mathcal{P}_{\lambda}(\widetilde{f}(\lambda,\cdot))(x).\;\Box$
	
	\indent Clearly the JEFT is the Poisson transform image of the Helgason Fourier transform and this relationship is the (noncompact) symmetric-space version of the earlier-noted classical case of the (full) Fourier transform on $\mathbb{R}^{n}$ being the one-dimensional completion of the classical Radon transform.
	
	\indent An immediate conclusion from the relationship in Lemma $2$ is to say that the JEFT, $f^{\triangle}(\lambda,x),$ is distinct from the Helgason Fourier transform, $\widetilde{f}(\lambda,b),$ but that they are the same {\it up to the Poisson transform, $\mathcal{P}_{\lambda}.$} In other words, the JEFT is the Poisson transform completion of the Helgason Fourier transform. However, {\it is this extension necessary at all?}
	
	\indent An affirmative response to the last question is deduced from the fact that, since $\mathcal{P}_{\lambda}:C(B)\rightarrow\mathcal{E}_{\lambda}(X)$ then, $$Df^{\triangle}(\lambda,\cdot)=\Gamma(D)(i\lambda)f^{\triangle}(\lambda,\cdot),$$ for every $D\in \textbf{D}(X).$ This is an earlier mentioned property of the classical Fourier transform on $\mathbb{R}^{n}$ which was not generally provable for the Helgason Fourier transform. This does not just show that the JEFT is the Poisson transform image of the Helgason Fourier transform but also that it is the (Poisson transform) completion of the Helgason Fourier transform. One fortunate twist of events is that the JEFT, being a function of $x,$ we can compute its Helgason Fourier transform as follows.
	
	\indent {\bf Lemma 3} ($[5.],$ p. $12$) For $f\in C^{\infty}_{c}(X),$ we have that $$\widetilde{f^{\triangle}}(\lambda,b)=\widetilde{f}(\lambda,b)\widehat{\varphi_{\mu}}(\lambda),$$ $\lambda,\mu\in\mathfrak{a}^{*}_{\mathbb{C}},$ $b\in B.\;\Box$
	
	\indent An immediate implication of Lemma $2$ is that, if $C^{\infty}_{c}(X)^{\triangle}$ denote the image of $C^{\infty}_{c}(X)$ under the JEFT, then its Helgason Fourier transform $$C^{\infty}_{c}(X)^{\triangle}\rightarrow\widetilde{(C^{\infty}_{c}(X)^{\triangle})}$$ is a projection into $C^{\infty}_{c}(\mathfrak{a}_{\mathbb{C}}\times X)$ with image within $\mathcal{H}(\mathfrak{a}_{\mathbb{C}}\times B)_{\mathfrak{w}},$ (See Theorem $5.1,$ of $[3]$ p. $270$ for meaning of notation) whose members split into a product of a member of $\mathcal{H}(\mathfrak{a}_{\mathbb{C}}\times B)_{\mathfrak{w}}$ and a member of the set $$\{\widehat{\varphi_{\mu}}(\lambda):\;\lambda\in\mathfrak{a}_{\mathbb{C}}\}=:S(\mu).$$ The conclusion is as follows.\\
	
	\indent {\bf Lemma 4}. $\widetilde{(C^{\infty}_{c}(X)^{\triangle})}
	\subseteq\mathcal{H}(\mathfrak{a}_{\mathbb{C}}\times B)_{\mathfrak{w}}\cdot S(\mu).\;\Box$
	
	\indent Is this inclusion strict? For which noncompact symmetric space $X$ is the reverse inclusion of Lemma $4$ valid? We may also ask in general; {\it What is the JEFT of the JEFT?} An important consequence of seeking the JEFT of the JEFT is in the formation of the Taylor-like power series for every $f\in C^{\infty}_{c}(X),$ as follows. Let $$(j_{1})(b)=f^{\triangle}(\lambda,b):=f^{\triangle}(\lambda,x)_{|_{x=b}},$$ $$(j_{2})(b)=(f^{\triangle}(\lambda,b))^{\triangle}:=
	(f^{\triangle}(\lambda,x))^{\triangle}_{|_{x=b}},\cdots,$$ with $(j_{o}f)(b):=f(b),$ in which $b\in B.$ The  formal power series $$\sum^{\infty}_{n=0}\frac{(j_{n})(b)}{n!}(x-b)^{n},\;x\in X,$$ may have something to say about the function $f,$ the symmetric space $X$ or the function-space $C^{\infty}_{c}(X).$
	
	\indent A preliminary study of the JEFT is contained in $[5.]$ and $[6.],$ including (in particular) its Plancherel formula and Paley-Wiener theorem. The Paley-Wiener theorem for the JEFT is worth quoting and it follows from that of Helgason Fourier transform and the injectivity of the Poisson transform. The result reads thus.
	
	\indent {\bf Theorem} ($[5.],$ p. $19$) Let $\lambda\in\mathfrak{a}^{*}_{\mathbb{C}}$ be simple. The Joint-Eigenspace Fourier transform on $X$ is a bijection of $C^{\infty}_{c}(X)$ onto the Hilbert space $H^{\infty}_{\lambda}(X).$ Moreover, we would have that $\psi(x)=f^{\triangle}(\lambda,x)$ is in $H^{\infty}_{\lambda}(X)$ iff $supp(f)\subset Cl(B_{R}(0)).\;\Box$
	
	\indent Hence, we now attain the commutative diagram
	
	 $$C^{\infty}_{c}(X)\rightarrow^{\widehat{f}(\lambda,b)}\rightarrow \mathcal{H}(\mathfrak{a}^{*}\times B)_{W}\rightarrow^{P_{\lambda}}\rightarrow \mathcal{E}_{\lambda}(X).$$
	 
	 \indent More fundamental questions about the JEFT have been considered in $[5.]$ and $[6.],$ through which we are now certain that {\it the Fourier transform on a noncompact symmetric space is the Joint-Eigenspace Fourier transform, the JEFT.} It would then mean that more conclusive results (than already known) are bound to emanate from the use of the JEFT in the harmonic analysis and representation theory of noncompact symmetric spaces.
	\ \\
	{\flushright{{\it Dedicated to the memory of Sigurdur Helgason}}}\\
	\ \\
	\ \\
	{\bf   References.}
	\begin{description}
		
		\item [{[1.]}] Barker, W. H.,  $L^p$ harmonic analysis on $SL(2, \R),$ \textit{Memoirs of American Mathematical Society,} $vol.$ {\bf 76}, no.: {\bf 393}. $1988$
		
		\item [{[2.]}] Camporesi, R., The Helgason Fourier transform for homogeneous vector bundles over compact Riemannian symmetric spaces - the local theory, {\it J. Funct. Anal.} {\bf 220} $2005,$ $97 - 117.$
		
		\item [{[3.]}] Helgason, S., \textit{Geometric analysis on symmetric spaces,} Mathematical Surveys and Monographs, vol. $39,$ Providence, Rhode Island $(1994)$
		
		\item [{[4.]}] Moharty, P., Ray, S. K., Sarkar, R. P., Sitaram, A., The Helgason-Fourier transform for symmetric spaces, $II,$ {\it J. Lie Theory,} {\bf 14}, $227 - 242.$
		
		\item [{[5.]}]  Oyadare, O. O.,  A Paley-Wiener theorem for the Joint-Eigenspace Fourier transform on noncompact symmetric spaces, {\it arXiv:$2409.09036.$} [math.FA], $(2024)$
		
		\item [{[6.]}]  Oyadare, O. O.,  A note on the $L^{2}-$harmonic analysis for the Joint-Eigenspace Fourier transform, {\it arXiv:$2410.09075.$} [math.FA], $(2024)$
		
		\item [{[7.]}] Sarkar, R. P., Sitaram, A., The Helgason Fourier transform for semisimple, Lie groups $I:$ the case of $SL(2,\mathbb{R}),$ {\it Bull. Austral. Math. Soc.,} {\bf 73}, $2006,$ $413 - 432.$
		
		\item [{[8.]}] Trombi, P. C. and Varadarajan, V. S., Spherical transforms on semisimple
		Lie groups, {\it Ann. Math.,} vol. {\bf 94.} $(1971),$ p. $246-303.$
		
		\item [{[9.]}] Varadarajan, V. S., {\it An introduction to harmonic analysis on semisimple
			Lie groups,} Cambridge Studies in Advanced Mathematics, {\bf 161,} Cambridge University Press, $1989.$
	
	\end{description}
\end{document}